\documentclass[letterpaper,11pt]{article}
\usepackage{amsmath, amsthm, amssymb}    % May not all be necessary
\usepackage{bridges}			% Custom bridges proceedings style
\usepackage{graphicx}			% For including pictures
\usepackage[colorlinks=true, urlcolor=blue, citecolor=black, linkcolor=black]{hyperref}  
\usepackage{subcaption}			% For captioning multi-panel figures

\usepackage{xfrac}

\newcommand{\Z}{\mathbb{Z}}
\newcommand{\R}{\mathbb{R}}
\newcommand{\C}{\mathbb{C}}

\newcommand{\SL}{\operatorname{SL}}
\newcommand{\diag}{\operatorname{diag}}
\newcommand{\figheight}{1.7in}

\title{Looping Animations Using the Modular Flow and Elliptic Functions}

\author{Clayton Shonkwiler
\vspace{10pt}\\
Department of Mathematics, Colorado State University; clayton.shonkwiler@colostate.edu} % end \author
\date{}					% Suppress any date on submissions

\usepackage{animate}

\begin{document}

\maketitle

% Prevent page number 1 from being printed on the first page.
\thispagestyle{empty}

\begin{abstract}

This paper describes an approach to generating looping animations using the modular flow and elliptic functions. The modular flow is a flow on lattices with many periodic orbits, and elliptic functions are meromorphic, doubly-periodic functions which can be visualized using domain coloring.

\end{abstract}

% Bridges papers are usually no more than 8 pages in length.  So
% there's really no need to have numbered sections, unless the
% author really needs to refer to sections by number within the paper's text.  
% So to suppress sequential section numbers, append an asterisk to 
% the \section command, as in:

%%%%%%%%%%%%%%%%%%%%%%%%%%%%%%%%%%%%%%%%%%%

I enjoy making looping animations using interesting mathematics and am always on the lookout for periodic orbits in spaces parametrizing mathematical objects. Every time I teach Complex Analysis I am inspired to try to make animations depicting interesting holomorphic or meromorphic functions, though I have often found it a challenge to come up with looping families of functions which are visually compelling. 

The goal here is to describe how to use the modular flow on the space of lattices to get periodic families of elliptic functions. Elliptic functions are meromorphic, doubly-periodic functions, so they are naturally associated with a planar lattice $\Lambda$ generated by the fundamental periods. In turn, the modular flow is a flow on the space of lattices which has many periodic orbits; these periodic orbits will give us our looping families of elliptic functions, which can then be visualized using domain coloring.\footnote{All figures in this paper are animations, which should be viewable in Acrobat Reader, Okular, and a few other PDF viewers. Other viewers will only show a single static frame. All animations can also be seen online at \url{https://shonkwiler.org/bridges26}. Source code for many of the animations in this paper can be found in the Geometry in Motion repository~\cite{GIM}.}

\section*{The Modular Flow}

% And don't try to mix and match numbered and unnumbered sections; it's one or the other.

I was first introduced to this topic by an animation of Etienne Jacob~\cite{EtienneJacob}. I will try to give a self-contained summary, but see~\cite{Bartlett,Ghys} for nice introductions. 

Consider the action of the non-zero real numbers $\R^*$ on the plane $\R^2$ by $s \cdot (x,y) = (s x, s^{-1} y)$. Typical orbits are hyperbolae of the form $xy = c$ for $c \neq 0$. %, though there are three special orbits, namely the $x$-axis minus the origin, the $y$-axis minus the origin, and the origin itself. % Extended to an action of $\C^*$ on $\CP^2$, this provides an elementary example motivating the construction of the Geometric Invariant Theory quotient~\cite{mumfordGeometricInvariantTheory1994}.
Since $\R^*$ is disconnected, it is convenient to modify this action by restricting to positive $s$ and thinking of $s = e^t$: then $t \cdot (x,y) = (e^t x, e^{-t} y)$ defines an action of $\R$ on $\R^2$. To get an action on lattices, recall that a lattice $\Lambda \subset \R^2$ consists of all points of the form $\{a\vec{u}_1 + b\vec{u}_2 : a,b \in \Z\}$ where $\vec{u}_1, \vec{u}_2 \in \R^2$ are linearly independent vectors. Given such a lattice $\Lambda$, we get a new lattice $t\cdot \Lambda$ for each $t \in \R$; this is the \emph{modular flow} on the space of lattices. 

Since the matrix $\diag(e^t,e^{-t})$ has determinant 1, its action is area-preserving, so the area of the fundamental domain of $t\cdot \Lambda$ is independent of $t$. In particular, we can interpret the generators $\vec{u}_1, \vec{u}_2$ of any unit-area $\Lambda$ as the columns of some $2 \times 2$ matrix $U$ with determinant 1, and then the generators of $t \cdot \Lambda$ are columns of $\diag(e^t,e^{-t}) U $, which still has determinant 1.

This action on matrices is very far from periodic, but it turns out that the action on lattices has some periodic orbits: for certain $\Lambda$ there exists $t_0>0$ so that $\Lambda = t_0 \cdot \Lambda$. As written, it's not so clear how to solve for $\Lambda$ and $t_0$; after all, if $\Lambda = t_0 \cdot \Lambda$, then for any vertex $\vec{v}$ of $\Lambda$ we have that $t_0 \cdot \vec{v}$ is also a vertex of $\Lambda$. But which one?

Notice that, if $B = \begin{bmatrix} a & b \\ c & d \end{bmatrix}$, then the columns of $UB$ will be $a \vec{u}_1 + c\vec{u}_2$ and $b\vec{u}_1 + d \vec{u}_2$. If $B$ is an integer matrix, then these will be elements of $\Lambda$, and if $B$ has determinant 1, then they will span a unit-area parallelogram. 

In other words, the group $\SL_2(\Z)$ of integer $2 \times 2$ matrices with determinant 1 sends the generators $\vec{u}_1$ and $\vec{u}_2$ of $\Lambda$ to different generators for $\Lambda$, so we can search for $\vec{u}_1$ and $\vec{u}_2$ so that $t_0 \cdot \vec{u}_1$ and $t_0 \cdot \vec{u}_2$ are the columns of $UB$ for some $B$. Equivalently, we can try to solve the matrix equation
\begin{equation}\label{eq:basic equation}
	UB = \diag(e^{t_0},e^{-t_0}) U
\end{equation}
for $U$ and $t_0$, where $B \in \SL_2(\Z)$ is chosen in advance.

For example, with $B = \begin{bmatrix} 2 & 1 \\ 1 & 1 \end{bmatrix}$ it's straightforward to check that, when $\varphi = \frac{1 + \sqrt{5}}{2}$ is the golden ratio, $t_0 = \ln(1+\varphi)$ and $U = \begin{bmatrix} \varphi & 1 \\ 1 & -\varphi \end{bmatrix}$ is a solution to \eqref{eq:basic equation}; see Figure~\ref{fig:lattice1}.

\begin{figure}[htbp]
	\centering
	\begin{minipage}[b]{\figheight} 
		\animategraphics[autoplay,loop,height=\figheight]{25}{lattice5/frame-}{1}{50}
	        	\subcaption{} % Add subcaption text if desired, or use \subcaption* to suppress (a), (b), etc. labels
	        	\label{fig:lattice1}
	\end{minipage}
	\qquad \qquad
	\begin{minipage}[b]{\figheight} 
		\animategraphics[autoplay,loop,height=\figheight]{25}{lattice6/frame-}{1}{50}
	        	\subcaption{} % Add subcaption text if desired, or use \subcaption* to suppress (a), (b), etc. labels
	        	\label{fig:lattice2}
	\end{minipage}
	\caption{Two periodic orbits in the space of lattices: (a) the lattice generated by $(\varphi, 1)$ and $(1,-\varphi)$, with period $\ln(1+\varphi) \approx 0.962$; (b) the lattice generated by $\left(2 \sqrt{\sfrac{2}{3}},-2 \sqrt{\sfrac{2}{3}}\right)$ and $(1,1)$, with period $\ln\left(5 + 2 \sqrt{6}\right)\approx 2.292$. Both animations take 2 seconds to show 1 period.}
	\label{fig:lattices}
\end{figure}

More generally, taking the transpose of \eqref{eq:basic equation} reveals that the rows of $U$ must be eigenvectors of $B^T$ with corresponding eigenvalues $e^{t_0}$ and $e^{-{t_0}}$. In other words, if $\lambda_1$ and $\lambda_2 = \frac{1}{\lambda_1}$ are the eigenvalues of $B$, then $t_0 = \pm \ln(\lambda_1)$ and the rows of $U$ are the corresponding eigenvectors. This explains golden ratio's appearance above: the characteristic polynomial of $\begin{bmatrix} 2 & 1 \\ 1 & 1 \end{bmatrix}$ is $x^2-3x+1$, which has roots $\lambda_1 = \varphi + 1$ and $\lambda_2 = 2-\varphi$.

The periodic flows starting from different lattices are qualitatively similar, but can already be reasonably visually compelling, as in the animation \emph{Modular Flow}; see Figure~\ref{fig:modularflow}.

\begin{figure}[htbp]
	\centering
		\animategraphics[autoplay,loop,height=\figheight]{25}{flow4/frame-}{1}{35}
	\caption{The animation \emph{Modular Flow}.}
	\label{fig:modularflow}
\end{figure}

\section*{Elliptic Functions}

While visualizing the modular flow on lattices is already interesting, we can take things a step further by introducing elliptic functions. To motivate this, identify $\R^2$ with the complex plane $\C$. We then consider doubly-periodic meromorphic functions\footnote{A \emph{meromorphic function} is complex-differentiable except at a discrete set of points at which it is allowed to blow up.} on $\C$. More precisely, if $\Lambda$ has fundamental domain with vertices $0, \omega_1, \omega_2$, and $\omega_1+\omega_2$ (here, $\omega_1, \omega_2 \in \C$ correspond to the basic vectors $\vec{u}_1, \vec{u}_2 \in \R^2$ defining $\Lambda$), then we are considering meromorphic functions $f$ on $\C$ so that
\[
	f(z + \omega_1) = f(z) = f(z+ \omega_2)
\]
for all $z \in \C$. Doubly-periodic meromorphic functions are usually called \emph{elliptic functions}; see Chapter 1 of~\cite{apostolModularFunctionsDirichlet1990} for a nice introduction and~\cite{lawdenEllipticFunctionsApplications1989} for a more comprehensive treatment.

In some sense the most basic elliptic function is the Weierstrass $\wp$-function; to emphasize the dependence on the lattice we will write this as $\wp_\Lambda$. The defining feature of $\wp_\Lambda$ is that its only poles are double poles at each lattice point and that the leading term of its Laurent series at 0 is $\frac{1}{z^2}$; here's the definition:
\[
	\wp_\Lambda(z) := \frac{1}{z^2} + \sum_{\omega \in \Lambda \backslash\{0\}}\left(\frac{1}{(z-\omega)^2} - \frac{1}{\omega^2}\right).
\]

In fact, every elliptic function $f$ with period lattice $\Lambda$ can be built out of $\wp_\Lambda$ and its derivative $\wp'_\Lambda$. More precisely, every such $f$ can be written as
\begin{equation}\label{eq:elliptic function rep}
	f(z) = R_1(\wp_\Lambda(z)) + \wp_\Lambda'(z)R_2(\wp_\Lambda(z)),
\end{equation}
where $R_1$ and $R_2$ are rational functions. 

Now, choose one of the lattices $\Lambda$ which has a periodic orbit under the modular flow and define an action of $\R$ on the Weierstrass function by $t \cdot \wp_\Lambda := \wp_{t \cdot \Lambda}$. Using~\eqref{eq:elliptic function rep}, we get an induced action on any elliptic function we like which will always have periodic orbits. 

\begin{figure}[htbp]
	\centering
	\begin{minipage}[b]{\figheight} 
		\animategraphics[autoplay,loop,height=\figheight]{25}{phase2/frame-}{1}{50}
	        	\subcaption{} % Add subcaption text if desired, or use \subcaption* to suppress (a), (b), etc. labels
	        	\label{fig:phase}
	\end{minipage}
	\qquad
	\begin{minipage}[b]{\figheight} 
		\animategraphics[autoplay,loop,height=\figheight]{25}{pprime/frame-}{1}{38}
	        	\subcaption{} % Add subcaption text if desired, or use \subcaption* to suppress (a), (b), etc. labels
	        	\label{fig:pprime}
	\end{minipage}
	\qquad
	\begin{minipage}[b]{\figheight} 
		\animategraphics[autoplay,loop,height=\figheight]{25}{jacobi/frame-}{1}{75}
	        	\subcaption{} % Add subcaption text if desired, or use \subcaption* to suppress (a), (b), etc. labels
	        	\label{fig:jacobi}
	\end{minipage}
	\caption{Three examples of periodic families of elliptic functions: (a) the animation \emph{Phase}, using the functions $\wp_{t \cdot \Lambda}$; (b) the functions $\wp'_{t \cdot \Lambda}$; (c) a family of Jacobi elliptic functions of type $\operatorname{cn}$ (see, e.g., \href{https://dlmf.nist.gov/22}{Chapter 22} of~\cite{NIST:DLMF}).}
	\label{fig:examples}
\end{figure}

With all this machinery in place, we can use domain coloring~\cite{farrisReviewsVisualComplex1998,farrisDomainColoringArgument2017} to visualize the periodic family of elliptic functions. Three examples are shown in Figure~\ref{fig:examples}, though a lot of experimentation is possible here. Three obvious parameters that can be modified are: (i) the choice of domain coloring palette; (ii) the speed of the animation; and (iii) the size of the visual field being displayed. For (i) I tend to use perceptually-uniform cyclic palettes, such as those found in~\cite{crameri_2023_8409685} and~\cite{CMasher}. For me, the main consideration in my choices for (ii) and (iii) is the balance between static (slow/zoomed-in) and dynamic (fast/zoomed-out), though file size restrictions are also an important constraint.

% \begin{figure}[htbp]
% 	\centering
% 		\animategraphics[autoplay,loop,height=2in]{25}{phase2/frame-}{1}{50}
% 	\caption{The animation \emph{Phase}~\cite{phase}.}
% 	\label{fig:phase}
% \end{figure}

In fact, this approach is not strictly limited to elliptic functions: we can use the same idea to get periodic families of quasi-periodic functions such as the Weierstrass $\sigma_\Lambda$ function, which is implicitly defined by the equation $\frac{d^2}{dz^2} \ln(\sigma_{\Lambda}) = -\wp_\Lambda$. See Figure~\ref{fig:sigma}.

\begin{figure}[htbp]
	\centering
		\animategraphics[autoplay,loop,height=\figheight]{25}{sigma/frame-}{1}{35}
	\caption{Domain coloring animation of the Weierstrass functions $\sigma_{t \cdot \Lambda}$.}
	\label{fig:sigma}
\end{figure}

% If $\Lambda$ is an honest lattice, then $\tau := \frac{\omega_2}{\omega_1} \notin \R$, and after changing the sign of $\omega_2$ if necessary, we can assume $\operatorname{Im}(\tau) > 0$. Indeed, we get an equivalent torus if we replace the pair $(\omega_1, \omega_2)$ with $(1,\tau)$

%%%%%%%%%%%%%%%%%%%%%%%%%%%%%%%%%%%%%%%
\section*{Acknowledgements}

Thanks to my Complex Analysis students and Anne Harding for motivation, and to the anonymous referees for helpful comments. This work was partially supported by the National Science Foundation (DMS–2107700). 

%%%%%%%%%%%%%%%%%%%%%%%%%%%%%%%%%%%%%%%
% References %
    
{\setlength{\baselineskip}{13pt} % tighten line spacing for bibliography
\raggedright				% no right justification for References
\bibliographystyle{bridges}
\bibliography{special,papers}
} % end setlength, raggedright
   
\end{document}